\documentclass[aps,prl,twocolumn,superscriptaddress,groupedaddress]{revtex4}

\usepackage[caption=false]{subfig}
\usepackage{graphicx} 
\usepackage{mathrsfs}
\usepackage{amsfonts}
\usepackage{amsmath}
\usepackage[noabbrev]{cleveref}
\usepackage[toc,page]{appendix}
\usepackage{import}
\usepackage[nice]{units}
\usepackage{float}
\usepackage{xcolor}
\usepackage{natbib}

\begin{document}
\title{Explicit, time-reversible and symplectic integrator for Hamiltonians in isotropic uniformly curved geometries}
\date{\today}
\author{Ana Silva}
\author{Eitan Ben Av}
\author{Efi Efrati} 
\email{efi.efrati@weizmann.ac.il} 
\affiliation{Department of
Physics of Complex Systems, Weizmann Institute of Science, Rehovot
76100, Israel}

\begin{abstract} 
The kinetic term of the $N$-body Hamiltonian system defined on the surface of the sphere is non-separable. As a result, standard explicit symplectic integrators are inapplicable. We exploit an underlying hierarchy in the structure of the kinetic term to construct an explicit time-reversible symplectic scheme of second order. We use iterative applications of the method to construct a fourth order scheme and demonstrate its efficiency. 
\end{abstract}
\maketitle

\section{\label{intro}Introduction\protect\\}
Most Hamiltonian systems cannot be solved analytically and their long time dynamics are thus commonly studied using numerical schemes that approximate the solutions of Hamilton's equations. Symplectic methods conserve the symplectic structure of the evolution equations and are therefore considered particularly suited and stable for numerically integrating autonomous Hamiltonian systems \cite{yoshida2,celestialsymplectic,channell1990symplectic,reich1996symplectic}. Simulating the long time dynamics of such systems requires not only stability and accuracy, but also efficiency, thus favoring explicit symplectic methods. 

Explicit symplectic methods are abundant for separable Hamiltonians where the kinetic term depends solely on the momenta and the potential term depends solely on the coordinates \citep{yoshida2,yoshida,forest1990fourth,Simhdynamics,GeometricIntBook,yoshida3,preto1999class}. Additional explicit symplectic schemes were devised for non-separable Hamiltonians of specific forms \cite{chin2009explicit,tao2016explicit2,tao2016explicit,benettin}. However, there is no explicit symplectic scheme for general Hamiltonians. The Hamiltonian of an autonomous $N$-body system in uniformly curved space, such as the surface of a sphere, is non-separable, and is not amenable to known existing explicit schemes. The obstruction to separability arises from the non-trivial metric coefficients in the kinetic term. We identify a hierarchy in the structure of the metric of uniform and isotropic spaces that allows the construction of an explicit, second order time reversible symplectic scheme. We use iterative applications of this method to obtain a fourth order integrator and evaluate its performance in comparison to other methods. 




\section{Constructing an explicit symplectic and time-reversible method}
The goal of the present section is to construct an explicit, symplectic, and time-reversible scheme  that approximates the exact flow map of a Hamiltonian with the following structure:
\begin{equation}
H(\mathbf{z})= H_1(p_{\theta}) + H_2(\theta,p_{\varphi}) + H_3(\theta,\varphi) \; ,
\end{equation}
with $\mathbf{z}=(\mathbf{q},\mathbf{p})^{T}=(\theta,\varphi,p_{\theta},p_{\varphi})^{T}$. 
The Hamiltonian $H(\mathbf{z})$ admits the splitting into the sum of several terms, but it is nonetheless non-separable. 
The non-separability arises due to the kinetic term $H_2(\theta,p_{\varphi})$, which prevents the complete splitting between terms that only depend on the position coordinates, and those that solely depend on the momentum coordinates. Such kinetic terms are typical for Hamiltonians defined on isotropic and uniform curved geometries (see SM). While the lack of full separability prevents the na\"ive implementation of standard splitting techniques, e.g. \cite{yoshida}, the hierarchical structure of the kinetic terms guides the design of a numerical scheme with the desired features along similar lines. To see this, let us start by examining Hamilton's equations:
\begin{equation}
\frac{d\mathbf{z}(t)}{dt} = \{ \mathbf{z}, H(\mathbf{z})\} = \hat{D}_H \mathbf{z},
\label{Hamilton_eq}
\end{equation}
where the curly braces represent the Poisson bracket, given by

\begin{equation}
\{ \mathbf{z}, H(\mathbf{z})\} =  \sum^2_{l=1} \frac{\partial \mathbf{z}}{\partial q_l}\frac{\partial H(\mathbf{z})}{\partial p_l} - \frac{\partial \mathbf{z}}{\partial p_l}\frac{\partial H(\mathbf{z})}{\partial q_l} \; .
\label{PoissonB}
\end{equation}
In eq.(\ref{Hamilton_eq}), the differential operator $\hat{D}_H$ is defined as $\hat{D}_H f=\{f, H(\mathbf{z})\}$, with $f$ an arbitrary function. Rewriting Hamilton's equations in operator form allows the introduction of a formal solution to the equations of motion~\cite{yoshida,SymplecticIntegratorsIntro}, namely

\begin{equation}
\mathbf{z}(\tau)= e^{\tau \hat{D}_H} \mathbf{z}(0) \; .
\label{formal_sol}
\end{equation}

From a computational stand point, the solution in eq. \eqref{formal_sol} is often impractical. Its implementation requires knowing the form of the flow map $e^{\tau \hat{D}_H}$, and therefore, requires knowing exactly how to solve the equations of motion. It, nevertheless, motivates the use of Yoshida's method~\cite{yoshida} in the construction of an approximate scheme for the flow map $e^{\tau \hat{D}_H}$. Recall that $H$ is given by the sum of three terms. As a consequence, the operator $\hat{D}_H$ also allows for a similar splitting, $\hat{D}_H=\hat{D}_{H_1}+\hat{D}_{H_2}+\hat{D}_{H_3}$ resulting in $e^{\tau \hat{D}_H}=e^{\tau (\hat{D}_{H_1}+\hat{D}_{H_2}+\hat{D}_{H_3})}$. We can thus aim at approximating the exact flow map $e^{\tau \hat{D}_H}$ as the successive concatenation of the individual flow maps $e^{\tau \lambda_j \hat{D}_{H_i}}$, with the number of concatenations $k$ and the value of the set of coefficients $\lambda_j=\{c_j,d_j,f_j\}$ determining the order of the integrator, $n$:
\begin{equation}
\begin{split}
&e^{\tau (\hat{D}_{H_1}+\hat{D}_{H_2}+\hat{D}_{H_3})}=\\ &= \prod^k_{j=1} e^{\tau c_j \hat{D}_{H_1}}e^{\tau d_j \hat{D}_{H_2}}e^{\tau f_j \hat{D}_{H_3}} + \mathcal{O}(\tau^{n+1}) \; .
\end{split}
\label{Yoshida_method}
\end{equation}
Whenever the individual Hamiltonian terms $H_i$ are integrable, each flow map $e^{\tau \hat{D}_{H_j}}$ is known exactly, and the method above can automatically be interpreted as the composition of symplectic maps, giving rise to a final method that is also symplectic~\cite{benettin,Simhdynamics}. In what follows, we show how to use Yoshida's framework to build explicit, symplectic and time-reversible methods of second and forth order.


A second order method is constructed from eq.(\ref{Yoshida_method}) by setting $k=3$, $c_1=d_1=d_2=c_3=\frac{1}{2}$, $f_1=1$ and $c_2=f_2=d_3=f_3=0$,
\begin{equation}
\begin{split}
&e^{\tau (\hat{D}_{H_1}+\hat{D}_{H_2}+\hat{D}_{H_3})}=\\ &= e^{\frac{\tau}{2} \hat{D}_{H_1}}e^{\frac{\tau}{2} \hat{D}_{H_2}}e^{\tau \hat{D}_{H_3}}e^{\frac{\tau}{2} \hat{D}_{H_2}}e^{\frac{\tau}{2} \hat{D}_{H_1}} + \mathcal{O}(\tau^{3}) \; .
\end{split}
\label{Yoshida_method1}
\end{equation}

To better understand the origins of eq.(\ref{Yoshida_method1}), we label $\hat{X}=\frac{\tau}{2} \hat{D}_{H_2}$, $\hat{Y}=\frac{\tau}{2} \hat{D}_{H_1}$ and $\hat{U}=\tau \hat{D}_{H_3}$, which allows the right hand side of eq.(\ref{Yoshida_method1}) to be expressed as the product $e^{\hat{Y}}e^{\hat{X}}e^{\hat{U}}e^{\hat{X}}e^{\hat{Y}}$. This product can be further simplified to $e^{\hat{Z'}}e^{\hat{T}}$, by setting $e^{\hat{Z'}}=e^{\hat{Y}}e^{\hat{X}}$ and $e^{\hat{T}}=e^{\hat{U}}e^{\hat{Z}}=e^{\hat{U}}e^{\hat{X}}e^{\hat{Y}}$. We are then left with solving the equation $e^{\hat{Z'}}e^{\hat{T}}=e^{\hat{W}}$, which immediately suggests the use of the Baker-Campbell-Hausdorff (BCH) formula to obtain a solution for $\hat{W}$. Given that the operators $\hat{X}$, $\hat{Y}$ and $\hat{U}$ are proportional to $\tau$, we only need to take the terms in the BCH formula that involve at most the product of two operators. Therefore, up to the desired accuracy, we get the following solution for $\hat{W}$:

\begin{equation}
\begin{split}
\hat{W}&=\hat{Z'}+\hat{T}+ \frac{1}{2}[\hat{Z'},\hat{T}]+ \mathcal{O}(\tau^{3})\\ &= 2\hat{X}+2\hat{Y}+\hat{U}
+ \mathcal{O}(\tau^{3}) \; .
\end{split}
\end{equation}

All the commutators arising from the term $[\hat{Z'},\hat{T}]$ cancel each other exactly, and the next contributions stemming from the BCH formula are already of the order $\tau^3$. We, thus, arrive at eq.(\ref{Yoshida_method1}). For this result to be useful, nonetheless, we need to be able to write down the individual flow maps $e^{\tau\hat{D}_{H_i}}$. We consider these to be given by the elementary one-step Euler method~\cite{Simhdynamics}:

\begin{equation}
e^{\tau \hat{D}_{H_i}}\mathbf{z}^{(j)}=\mathbf{z}^{(j+1)}=\mathbf{z}^{(j)} + \tau \hat{D}_{H_i} \mathbf{z}^{(j)} \; ,
\end{equation}

with the action of the operator $\hat{D}_{H_i}$ as given in eq.(\ref{Hamilton_eq}). By acting successively with the individual flow maps $e^{\tau\hat{D}_{H_i}}$ on $\mathbf{z}^{(j)}$, as described in eq.(\ref{Yoshida_method1}), we obtain the following integrator:

\begin{equation}
\begin{split}
&\theta^{(j+\frac{1}{2})}=\theta^{(j)} + \frac{\tau}{2} \frac{\partial H_1}{\partial p_{\theta}} \big(p^{(j)}_{\theta}\big)\\
&\varphi^{(j+\frac{1}{2})} = \varphi^{(j)} + \frac{\tau}{2} \frac{\partial H_2}{\partial p_{\varphi}} \big(\theta^{(j+\frac{1}{2})},p^{(j)}_{\varphi} \big) \\
& p^{(j+1)}_{\varphi}= p^{(j)}_{\varphi} - \tau \frac{\partial H_3}{\partial \varphi} \big(\theta^{(j+\frac{1}{2})},\varphi^{(j+\frac{1}{2})} \big) \\
& p^{(j+1)}_{\theta}= p^{(j)}_{\theta} - \tau \Bigg( \frac{\partial H_3}{\partial \theta} \big(\theta^{(j+\frac{1}{2})},\varphi^{(j+\frac{1}{2})} \big) + \\ & + \frac{1}{2} \frac{\partial H_2}{\partial \theta} \big(\theta^{(j+\frac{1}{2})},p^{(j)}_{\varphi} \big) + \frac{1}{2} \frac{\partial H_2}{\partial \theta} \big(\theta^{(j+\frac{1}{2})},p^{(j+1)}_{\varphi} \big)  \Bigg) \\
&\theta^{(j+1)}=\theta^{(j+\frac{1}{2})} + \frac{\tau}{2} \frac{\partial H_1}{\partial p_{\theta}} \big(p^{(j+1)}_{\theta}\big) \\
&\varphi^{(j+1)} = \varphi^{(j+\frac{1}{2})} + \frac{\tau}{2} \frac{\partial H_2}{\partial p_{\varphi}} \big(\theta^{(j+\frac{1}{2})},p^{(j+1)}_{\varphi} \big) \; .
\end{split}
\label{final_method}
\end{equation}

It then follows that the proposed numerical integrator, when executed in the order of instructions listed in eq.(\ref{final_method}), is entirely composed from explicit steps. This greatly simplifies its numerical implementation.

We may write the integrator in eq.(\ref{final_method}) in a more compact form, namely as the set of instructions

\begin{equation}
\begin{split}
& \mathbf{z}^{(j+\frac{1}{2})}=\mathbf{z}^{(j)} + \frac{\tau}{2} J \nabla_{\mathbf{z}} H \big( \mathbf{q}^{(j+\frac{1}{2})},\mathbf{p}^{(j)} \big) \\
& \mathbf{z}^{(j+1)}= \mathbf{z}^{(n+\frac{1}{2})} + \frac{\tau}{2} J \nabla_{\mathbf{z}} H \big( \mathbf{q}^{(j+\frac{1}{2})},\mathbf{p}^{(j+1)} \big) \; ,
\end{split}
\label{final_method1}
\end{equation}
where $J$ is the canonical structure matrix, given by $J=\begin{pmatrix} \mathbf{0} & \mathbf{1} \\ -\mathbf{1} &\mathbf{0}  \end{pmatrix}$, with $\mathbf{0}$ representing a $2\times2$ matrix filled with zeros, and $\mathbf{1}$ the $2\times2$ identity matrix. Although the more compact notation used in eq.(\ref{final_method1}) hides the fact that integrator can be made explicit, it allows us to deduce the symplectic and time-reversible properties of the method. Written in the compact notation of eq.(\ref{final_method1}), the method proposed in eq.(\ref{Yoshida_method1}) can be identified with the generalised St\"{o}rmer-Verlet/leapfrog method, which is known to be a second order symplectic and time-reversible method~\cite{Simhdynamics,GeometricIntBook,hairer2003geometric}. Concatenating three such steps with appropriately chosen time differences, namely 
\[
\tilde{\tau}_1=\frac{\tau}{2-2^{1/3}},\quad
\tilde{\tau}_2=-\frac{2^{1/3}\tau}{2-2^{1/3}},\quad
\tilde{\tau}_3=\frac{\tau}{2-2^{1/3}},
\]
results in a fourth order integration scheme advancing the system by the time step $\tau$ \cite{yoshida}.

\section{Comparison of numerical code}
\begin{figure*}
    \includegraphics[width=\textwidth]{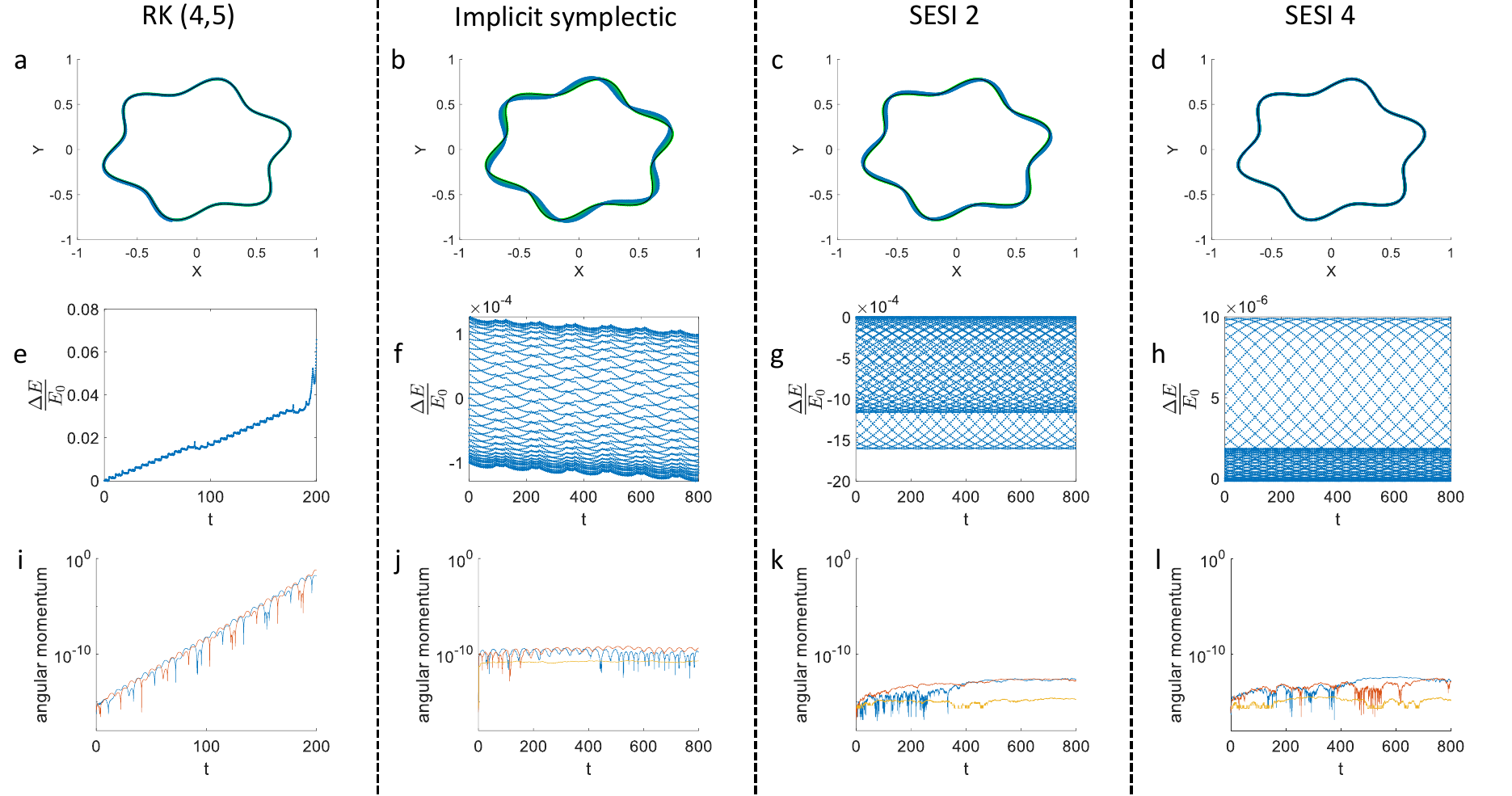}
    \caption[width=\textwidth]{\textbf{Comparison of the numerical integration schemes:} Simulating the time evolution of the Hamiltonian in equation \ref{eqn:Hamiltonian_sim}, for special symmetrical initial conditions. The exact solution, given in the SM\ref{SM}, displays a closed trajectory whose projection onto the XY plane yields a 6-petaled shape. 
    The schemes compared are:
    A commercially available explicit Runge-Kutta (4,5)\cite{dormand1980family,shampine1997matlab} (first column), an implicit symplectic midpoint second order scheme (second column) \cite{BLSR2004}, and the $2^{nd}$ and $4^{th}$ order symplectic explicit spherical integrator (SESI) introduced in this work (third and fourth columns, respectively). (a),(b),(c),(d) - 
    Depict the projection of the trajectory of one of the masses on the XY plane.  The results of the numerical integrators shown in (green and blue) dots overlaid by the exact solution, plotted in a continuous (black) line. 
    The shading (color) of the dots vary with time from dark (green) to light (blue). (e),(f),(g),(h) - Display energy deviation, where $\Delta E\equiv E-E_0$. (i),(j),(j),(k) - Present the numerical angular momentum drift on $x,y,z$ axes, in blue, red, and yellow respectively.}
    \label{fig:comparison_panel}
\end{figure*}

To benchmark our integrator, we simulate the following Hamiltonian 
\begin{equation}
\begin{split}
H	=\sum_{i}\frac{1}{2m}\left(p_{i,\theta}^{2}+\left(\frac{p_{i,\phi}}{\sin\left(\theta_i\right)}\right)^{2}\right)\\
	+\sum_{i,j}\left(\frac{\sqrt{1-\left(\frac{2\sin\left(\frac{L_{ij}}{2R}\right)}{\sqrt{3}}\right)^{2}}-\cos\left(\theta_{0}\right)}{\frac{2\sin\left(\frac{L_{ij}}{2R}\right)}{\sqrt{3}}}\right)^{2},
\end{split}
\label{eqn:Hamiltonian_sim}
\end{equation}
which allows an explicit solution (See SM\ref{SM} for details). By fine tuning the initial conditions we obtain a closed trajectory in real space, which enables an easier visualization of the deviations from the exact solution in long time trajectories. Figure \ref{fig:comparison_panel} shows the comparison of the Symplectic Explicit Spherical Integrator (SESI) derived here for both second order and fourth order schemes, with the implicit midpoint symplectic second order scheme \cite{BLSR2004} (see SM\ref{SM} for details), as well as with a commercially available explicit Runge-Kutta (RK) method based on the Dormand-Prince (4,5) pair \cite{dormand1980family,shampine1997matlab}.  We observe that both the $2^{nd}$ and the $4^{th}$ order SESI schemes demonstrate a bounded error in energy, while the non symplectic RK(4,5) displays an unbounded error in energy that grows with time. The error in angular momentum is unbounded, yet better in the symplectic integrators than the RK(4,5). Moreover, we observe that the $4^{th}$ order symplectic integrator conserves the closed shape of the solution better than the other solvers, for the given time window. The run time of the implicit symplectic solver was two orders of magnitude longer than the comparable $2^{nd}$ order SESI.

\section{discussion}
Efficient and reliable long time numerical integration of Hamiltonian systems favor explicit symplectic schemes. The kinetic term in many-body Hamiltonians in uniformly curved space is non-seperable, precluding the usage of existing explicit symplectic methods. We exploit the geometric hierarchy in the structure of the kinetic term on the surface of the sphere, to produce an explicit symplectic scheme. We note that while the resulting integrator is symplectic and energy deviations remain bounded, other conserved quantities arising from the symmetry of the metric, such as the angular momentum, are not identically conserved.  

These results may be generalized to higher dimensions as well as additional Hamiltonians with a certain structure.  The metrics of all uniformly curved spaces of $N$-dimensions share a similar hierarchy to the one we exploited above; namely that the kinetic term for each of the masses can be written as the sum of $N$ terms $K_i$, where each $K_i$ depends on $p_i$ and on the preceding coordinates $q_1,q_2,\dots, q_{i-1}$:
\[
K=\sum_{i=1}^N K_i(q_1,q_2,\dots, q_{i-1},p_i).
\]
Whenever the non-separable kinetic term displays this form, we can utilize a similar maneuver to the one employed in this work to produce an explicit symplectic scheme as expressed in equation \eqref{final_method1}: Advancing the $i^{th}$ coordinate, $q_i$, from time $t$ to time $t+\tau/2$ is performed using the knowledge of $q_{k<i}$  at time $t+\tau/2$, and of $p_i$ at time $t$.
Due to the hierarchical structure described above, this could be carried out explicitly.  Advancing the corresponding momenta to time $t+\tau/2$ can then be carried out explicitly as they too depend only on the momenta at time $t$ and on the already calculated coordinates at time $t+\tau/2$. To obtain time-reversibility we next perform the mirrored version of the above step: We note that the hierarchy above implies that $\partial H/\partial q_i$ may only depend on momenta terms $p_k$ such that $i<k$. Thus, starting with the momenta of largest index, we may successively  advance the momenta from time $t+\tau/2$ to time $t+\tau$ using the knowledge of the coordinates at time $t+\tau/2$ and the already calculated momenta at time $t+\tau$. We conclude by advancing the coordinates from time $t+\tau/2$ to time $t+\tau$ by using the coordinates at time $t+\tau/2$ and the already calculated momenta at time $t+\tau$. We thus obtain an explicit, time-reversible, symplectic, integration-scheme. The resulting second order scheme may be improved to yield higher even orders using standard procedures \cite{yoshida}. 


\bibliography{BibSymplecticIntegrators}

\begin{thebibliography}{19}
\expandafter\ifx\csname natexlab\endcsname\relax\def\natexlab#1{#1}\fi
\expandafter\ifx\csname bibnamefont\endcsname\relax
  \def\bibnamefont#1{#1}\fi
\expandafter\ifx\csname bibfnamefont\endcsname\relax
  \def\bibfnamefont#1{#1}\fi
\expandafter\ifx\csname citenamefont\endcsname\relax
  \def\citenamefont#1{#1}\fi
\expandafter\ifx\csname url\endcsname\relax
  \def\url#1{\texttt{#1}}\fi
\expandafter\ifx\csname urlprefix\endcsname\relax\def\urlprefix{URL }\fi
\providecommand{\bibinfo}[2]{#2}
\providecommand{\eprint}[2][]{\url{#2}}

\bibitem[{\citenamefont{Yoshida}(1993)}]{yoshida2}
\bibinfo{author}{\bibfnamefont{H.}~\bibnamefont{Yoshida}},
  \bibinfo{journal}{Qualitative and Quantitative Behaviour of Planetary
  Systems} pp. \bibinfo{pages}{27--43} (\bibinfo{year}{1993}).

\bibitem[{\citenamefont{Gladman et~al.}(1991)\citenamefont{Gladman, Duncan, and
  Candy}}]{celestialsymplectic}
\bibinfo{author}{\bibfnamefont{B.}~\bibnamefont{Gladman}},
  \bibinfo{author}{\bibfnamefont{M.}~\bibnamefont{Duncan}}, \bibnamefont{and}
  \bibinfo{author}{\bibfnamefont{J.}~\bibnamefont{Candy}},
  \bibinfo{journal}{Celestial Mechanics and Dynamical Astronomy}
  \textbf{\bibinfo{volume}{52}}, \bibinfo{pages}{221} (\bibinfo{year}{1991}).

\bibitem[{\citenamefont{Channell and Scovel}(1990)}]{channell1990symplectic}
\bibinfo{author}{\bibfnamefont{P.~J.} \bibnamefont{Channell}} \bibnamefont{and}
  \bibinfo{author}{\bibfnamefont{C.}~\bibnamefont{Scovel}},
  \bibinfo{journal}{Nonlinearity} \textbf{\bibinfo{volume}{3}},
  \bibinfo{pages}{231} (\bibinfo{year}{1990}).

\bibitem[{\citenamefont{Reich}(1996)}]{reich1996symplectic}
\bibinfo{author}{\bibfnamefont{S.}~\bibnamefont{Reich}}, \bibinfo{journal}{SIAM
  journal on numerical analysis} \textbf{\bibinfo{volume}{33}},
  \bibinfo{pages}{475} (\bibinfo{year}{1996}).

\bibitem[{\citenamefont{Haruo}(1990)}]{yoshida}
\bibinfo{author}{\bibfnamefont{Y.}~\bibnamefont{Haruo}},
  \bibinfo{journal}{Physics letters A} \textbf{\bibinfo{volume}{150}},
  \bibinfo{pages}{262} (\bibinfo{year}{1990}).

\bibitem[{\citenamefont{Forest and Ruth}(1990)}]{forest1990fourth}
\bibinfo{author}{\bibfnamefont{E.}~\bibnamefont{Forest}} \bibnamefont{and}
  \bibinfo{author}{\bibfnamefont{R.~D.} \bibnamefont{Ruth}},
  \bibinfo{journal}{Physica D: Nonlinear Phenomena}
  \textbf{\bibinfo{volume}{43}}, \bibinfo{pages}{105} (\bibinfo{year}{1990}).

\bibitem[{\citenamefont{Leimkuhler and
  Reich}(2004{\natexlab{a}})}]{Simhdynamics}
\bibinfo{author}{\bibfnamefont{B.}~\bibnamefont{Leimkuhler}} \bibnamefont{and}
  \bibinfo{author}{\bibfnamefont{S.}~\bibnamefont{Reich}},
  \emph{\bibinfo{title}{Simulating hamiltonian dynamics}}, \bibinfo{number}{14}
  (\bibinfo{publisher}{Cambridge university press},
  \bibinfo{year}{2004}{\natexlab{a}}).

\bibitem[{\citenamefont{Hairer et~al.}(2006)\citenamefont{Hairer, Lubich, and
  Wanner}}]{GeometricIntBook}
\bibinfo{author}{\bibfnamefont{E.}~\bibnamefont{Hairer}},
  \bibinfo{author}{\bibfnamefont{C.}~\bibnamefont{Lubich}}, \bibnamefont{and}
  \bibinfo{author}{\bibfnamefont{G.}~\bibnamefont{Wanner}},
  \emph{\bibinfo{title}{Geometric numerical integration: structure-preserving
  algorithms for ordinary differential equations}}, vol.~\bibinfo{volume}{31}
  (\bibinfo{publisher}{Springer Science \& Business Media},
  \bibinfo{year}{2006}).

\bibitem[{\citenamefont{Yoshida}(1992)}]{yoshida3}
\bibinfo{author}{\bibfnamefont{H.}~\bibnamefont{Yoshida}}, in
  \emph{\bibinfo{booktitle}{Symposium-International Astronomical Union}}
  (\bibinfo{organization}{Cambridge University Press}, \bibinfo{year}{1992}),
  vol. \bibinfo{volume}{152}, pp. \bibinfo{pages}{407--411}.

\bibitem[{\citenamefont{Preto and Tremaine}(1999)}]{preto1999class}
\bibinfo{author}{\bibfnamefont{M.}~\bibnamefont{Preto}} \bibnamefont{and}
  \bibinfo{author}{\bibfnamefont{S.}~\bibnamefont{Tremaine}},
  \bibinfo{journal}{The Astronomical Journal} \textbf{\bibinfo{volume}{118}},
  \bibinfo{pages}{2532} (\bibinfo{year}{1999}).

\bibitem[{\citenamefont{Chin}(2009)}]{chin2009explicit}
\bibinfo{author}{\bibfnamefont{S.~A.} \bibnamefont{Chin}},
  \bibinfo{journal}{Physical Review E} \textbf{\bibinfo{volume}{80}},
  \bibinfo{pages}{037701} (\bibinfo{year}{2009}).

\bibitem[{\citenamefont{Tao}(2016{\natexlab{a}})}]{tao2016explicit2}
\bibinfo{author}{\bibfnamefont{M.}~\bibnamefont{Tao}},
  \bibinfo{journal}{Journal of Computational Physics}
  \textbf{\bibinfo{volume}{327}}, \bibinfo{pages}{245}
  (\bibinfo{year}{2016}{\natexlab{a}}).

\bibitem[{\citenamefont{Tao}(2016{\natexlab{b}})}]{tao2016explicit}
\bibinfo{author}{\bibfnamefont{M.}~\bibnamefont{Tao}},
  \bibinfo{journal}{Physical Review E} \textbf{\bibinfo{volume}{94}},
  \bibinfo{pages}{043303} (\bibinfo{year}{2016}{\natexlab{b}}).

\bibitem[{\citenamefont{Benettin et~al.}(2001)\citenamefont{Benettin,
  Cherubini, and Fass{\`o}}}]{benettin}
\bibinfo{author}{\bibfnamefont{G.}~\bibnamefont{Benettin}},
  \bibinfo{author}{\bibfnamefont{A.~M.} \bibnamefont{Cherubini}},
  \bibnamefont{and}
  \bibinfo{author}{\bibfnamefont{F.}~\bibnamefont{Fass{\`o}}},
  \bibinfo{journal}{SIAM Journal on Scientific Computing}
  \textbf{\bibinfo{volume}{23}}, \bibinfo{pages}{1189} (\bibinfo{year}{2001}).

\bibitem[{\citenamefont{Donnelly and
  Rogers}(2005)}]{SymplecticIntegratorsIntro}
\bibinfo{author}{\bibfnamefont{D.}~\bibnamefont{Donnelly}} \bibnamefont{and}
  \bibinfo{author}{\bibfnamefont{E.}~\bibnamefont{Rogers}},
  \bibinfo{journal}{American Journal of Physics} \textbf{\bibinfo{volume}{73}},
  \bibinfo{pages}{938} (\bibinfo{year}{2005}).

\bibitem[{\citenamefont{Hairer et~al.}(2003)\citenamefont{Hairer, Lubich, and
  Wanner}}]{hairer2003geometric}
\bibinfo{author}{\bibfnamefont{E.}~\bibnamefont{Hairer}},
  \bibinfo{author}{\bibfnamefont{C.}~\bibnamefont{Lubich}}, \bibnamefont{and}
  \bibinfo{author}{\bibfnamefont{G.}~\bibnamefont{Wanner}},
  \bibinfo{journal}{Acta numerica} \textbf{\bibinfo{volume}{12}},
  \bibinfo{pages}{399} (\bibinfo{year}{2003}).

\bibitem[{\citenamefont{Dormand and Prince}(1980)}]{dormand1980family}
\bibinfo{author}{\bibfnamefont{J.~R.} \bibnamefont{Dormand}} \bibnamefont{and}
  \bibinfo{author}{\bibfnamefont{P.~J.} \bibnamefont{Prince}},
  \bibinfo{journal}{Journal of computational and applied mathematics}
  \textbf{\bibinfo{volume}{6}}, \bibinfo{pages}{19} (\bibinfo{year}{1980}).

\bibitem[{\citenamefont{Shampine and Reichelt}(1997)}]{shampine1997matlab}
\bibinfo{author}{\bibfnamefont{L.~F.} \bibnamefont{Shampine}} \bibnamefont{and}
  \bibinfo{author}{\bibfnamefont{M.~W.} \bibnamefont{Reichelt}},
  \bibinfo{journal}{SIAM journal on scientific computing}
  \textbf{\bibinfo{volume}{18}}, \bibinfo{pages}{1} (\bibinfo{year}{1997}).

\bibitem[{\citenamefont{Leimkuhler and Reich}(2004{\natexlab{b}})}]{BLSR2004}
\bibinfo{author}{\bibfnamefont{B.}~\bibnamefont{Leimkuhler}} \bibnamefont{and}
  \bibinfo{author}{\bibfnamefont{S.}~\bibnamefont{Reich}},
  \emph{\bibinfo{title}{Simulating hamiltonian dynamics}}, \bibinfo{number}{14}
  (\bibinfo{publisher}{Cambridge university press},
  \bibinfo{year}{2004}{\natexlab{b}}).

\end{thebibliography}

\newpage


\newpage
\onecolumngrid


\section{\label{SM} \large S\lowercase{upplementary} M\lowercase{aterial}: \protect\\
E\lowercase{xplicit, time-reversible and symplectic integrator for} H\lowercase{amiltonians in isotropic uniformly curved geometries}}


\section{Order of the integrator}
\begin{figure}[h]
    \includegraphics[width=0.5\linewidth]{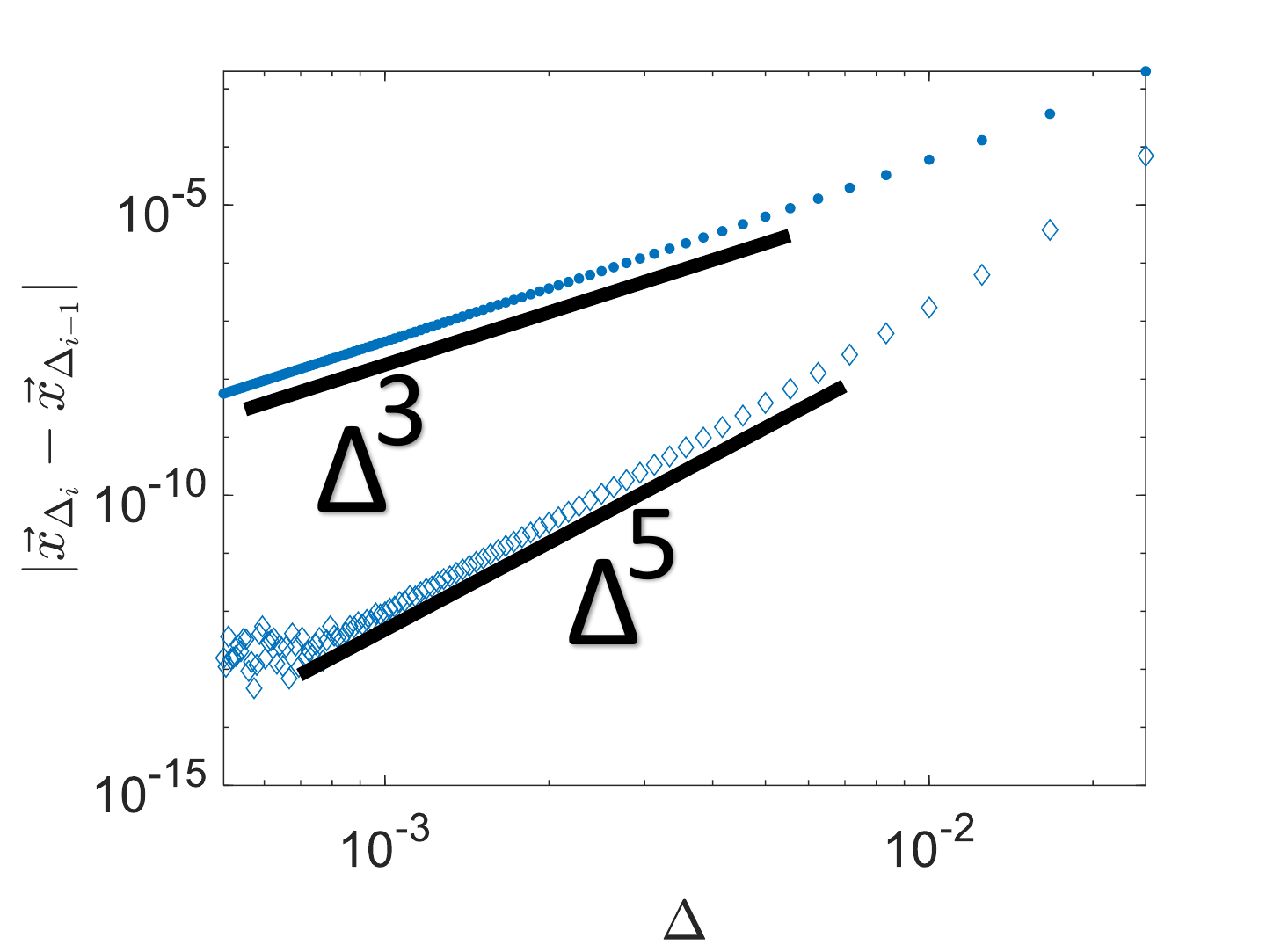}
    \caption[width=0.5\linewidth]{\textbf{Numerical schemes convergence rate.} $\vec{x}_{\Delta_{i}}\equiv \left(\theta(t=10),\phi(t=10),p_{\theta}(t=10),p_{\phi}(t=10)\right)$ calculated using a time step of $\Delta_{i}$. Filled dots correspond to SESI 2; hollow diamonds correspond to SESI 4.}
    \label{fig:conv}
\end{figure}

For the purpose of verifying the order of the integration scheme we run the simulation to a fixed final time $t=10$ using different time steps $\Delta t$, and examine the difference between the value of the final position $x(t=10)$ for the different time steps. Figure \ref{fig:conv} shows the difference in the final position $x(t=10)$ between two runs with nearby time steps, as a function of the time step magnitude indicating that the SESI 2 and SESI 4 schemes are indeed second and fourth order respectively.

\section{Finding an explicit solution}
To benchmark our numerical schemes' accuracy we seek to find a potential that will allow an explicit analytic solution. We seek this solution as a highly symmetric quasi-one-dimensional solution for a specially designed interaction.

We start by considering the three mass system on the sphere. Each of the masses position is given by 
\[
\mathbf{r}_i=(\sin(\theta_i)\cos(\phi_i),\sin(\theta_i) \sin(\phi_i),\cos(\theta_i)).
\]
We assume the sought solution is associated with the three fold symmetry for the problem, i.e.
\[
\theta_i=\theta,\quad \phi_1=\phi,\quad
\phi_2=\phi+2\pi/3,\quad
\phi_3=\phi+4\pi/3.
\]
The geodesic (or any other) distance between every two particles becomes in this highly symmetric case only a function of $\theta$. The real space distance simply measures $\sqrt{3}\sin(\theta)$. The kinetic terms are proportional to $\dot{\theta}^2+\sin(\theta)^2\dot\phi^2$. The Lagrangian of the system thus reduces (after dividing by the constant $3m/2R^2$):
\[
\mathcal{L}=\dot{\theta}^2+\sin(\theta)^2\dot\phi^2-V(\theta).
\]
The conserved angular momentum of the system is given by
\[
L=\sin(\theta)^2 \dot\phi.
\]
The Euler Lagrange equation thus yields
\[
2 \sin(\theta)^{-2} \cot(\theta)L^2-V'(\theta)-2 \ddot\theta 
=0,
\]
which could be multiplied by $\dot{\theta}$ and explicitly integrated to give
\[
V(\theta)+L^2 \cot(\theta)^{2}+\dot{\theta}^2=c_1.
\]
We may now identify the system's conserved energy and rewrite the equation as
\[
V(\theta)+L^2 (\sin(\theta)^{-2}-1)+\dot{\theta}^2=c_1=E_0-L^2.
\]
Setting $\psi=\cos(\theta)$ we have
\[
\dot{\psi}^2=E_0-L^2 -E_0\psi^2-U(\psi) 
\]
where 
$U(\psi)=U(\cos(\theta))=V(\theta)\sin(\theta)^{2}$.
We are thus required to integrate
\[
\int \frac{d\psi}{\sqrt{E_0-L^2 -E_0\psi^2-U(\psi)}}=t.
\]
What remains is to choose $V(\theta)$ and the corresponding $U(\psi)$ to render the equations solvable. One very convenient choice for the interactions reads 
\[
V(\theta)=(\cos(\theta)-\cos(\theta_0))^2/\sin(\theta)^2.
\]
While this potential diverges near the origin, it gives a nice parabolic behavior in the vicinity of $\theta_0$, which is where we make use of it. Moreover, this potential yields $U(\psi)=(\psi-\psi_0)^2$. It is now straightforward to explicitly calculate the full time evolution for $\psi$:
\[
\psi=\frac{\psi_0}{1+E_0}+\sqrt{\frac{E_0-L^2-\frac{E_0\psi_0^2}{1+E_0}}{1+E_0}}\sin(\sqrt{1+E_0} t).
\]
Writing the three dimensional distance between every two masses as
$L^{3D}=\sqrt{3}R\sin(\theta)$ and its relation to the corresponding geodesic distance on the sphere through 
\[
L^S=2R \arcsin(L^{3D}/2R),
\]
yields the potential in the form (11).

\section{Simulation details}
We set $\theta_0 = \frac{\pi}{4}$. In order to get a closed shape trajectory in real space, we choose $p_{\theta}=\frac{1}{4}$ and $p_{\phi} = 0.1727174029854043$ such that $\int_{0}^{\frac{2\pi}{\sqrt{E_{0}+1}}}\dot{\varphi}dt=\frac{2\pi}{6}$. The integration time step is $\tau = 0.1$ and the total integration time is $t=800=8000\tau$.

The implicit midpoint symplectic second order integration scheme is [19]
\begin{equation}
\begin{split}
    p_{n+\frac{1}{2}}	=p_{n}-\frac{\tau}{2}\frac{\partial H}{\partial q}\left(q_{n+\frac{1}{2}},p_{n+\frac{1}{2}}\right)\\
q_{n+\frac{1}{2}}	=q_{n}+\frac{\tau}{2}\frac{\partial H}{\partial p}\left(q_{n+\frac{1}{2}},p_{n+\frac{1}{2}}\right)\\
q_{n+1}	=q_{n+\frac{1}{2}}+\frac{\tau}{2}\frac{\partial H}{\partial p}\left(q_{n+\frac{1}{2}},p_{n+\frac{1}{2}}\right)\\
p_{n+1}	=p_{n+\frac{1}{2}}-\frac{\tau}{2}\frac{\partial H}{\partial q}\left(q_{n+\frac{1}{2}},p_{n+\frac{1}{2}}\right)
\end{split}
\end{equation}

\end{document}